# Classical Solutions for a Class of Burgers Equations

Svetlin G. Georgiev and Gal Davidi

Wed Dec 16 08:44:08 2020

In this paper we consider a class of Burgers equations. We propose a new method for investigation for existence of classical solutions.

## 1    Introduction

When solving a wide range of problems related to nonlinear acoustics, we may describe the nonlinear sound waves in fluids by using the Burgers model equation. This equation is named after Johannes Martinus Burgers, who published an equation of this form in his paper concerning turbulent phenomena modelling in 1948 (see [1]).However, this type of equation used in continuum mechanics was first presented in a paper published in ameteorology journal by H. Bateman in 1915 (see [2]). Thanks to the fact that this journal was not read by experts in continuum mechanics, the equation has become known as the Burgers equation. The possibilities of using this equation in nonlinear acoustics were probably discovered by J. Cole in 1949. Subsequently E. Hopf (in 1950) and J. Cole (in 1951) published the transformation independently in their papers (see [3, 4]). The transformation enables us to find the general analytic solution of the Burgers equation, as a result of which this equation plays a major role in nonlinear acoustics. The Burgers equation can be used for modelling both travelling and standing nonlinear plane waves [5-8]. The equation is the simplest model equation that can describe the second order nonlinear effects connected with the propagation of high-amplitude (finite-amplitude waves) plane waves and, in addition, the dissipative effects in real fluids [9-12]. There are several approximate solutions of the Burgers equation (see [13]). These solutions are always fixed to areas before and after the shock formation. For an area where the shock wave is forming no approximate solution has yet been found. It is therefore necessary to solve the Burgers equation numerically in this area (see [14,15]). Numerical solutions themselves have difficulties with stability and accuracy.

Here we focus our attention on a class of Burgers equations and we will investigate it for existence of classical solutions. More precisely, consider the IVP

$$u_t + uu_x = 0, \quad t > 0, \quad x \in \mathbb{R},$$

$$u(0, x) = \phi(x), \tag{1.1}$$

where

$$\phi \in \cup_{j=-\infty}^{\infty} \mathcal{C}_0^1([j, j+1]), \quad \parallel \phi \parallel_{L^2(\mathbb{R})} \leq 1,$$

$$\sup_{x \in \mathbb{R}} |\phi(x)|, \quad \sup_{x \in \mathbb{R}} |\phi'(x)| < 1, \qquad (1.2)$$

$$|\phi|, \quad |\phi'| \leq \frac{1}{2^{|j|+1}} \quad on \quad [j, j+1], \quad j \in \mathbb{Z}.$$

Our main result is as follows.

**Theorem 1.1** *Assume (1.2). Then the IVP (1.1) has a solution* $u \in \mathcal{C}^1(\mathbb{R}, \cup_{j=-\infty}^{\infty} \mathcal{C}_0^1([j, j+1]))$.

The paper is organized as follows. In the next section we will give some preliminaries which will be used for the proof of our main result. In Section 3 we give the proof of our main result. In Section 4 we make some conclusions.

## 2   Preliminaries

To prove our main result we will use a fixed point theorem for a sum two operators one of which is an expansive operator.

**Theorem 2.1 (, Theorem 2.4)**  *Let  $X$  be a nonempty closed convex subset of a Banach space  $E$ . Suppose that  $T$  and  $S$  map  $X$  into  $E$  such that*

1.  $S$  is continuous and  $S(X)$  resides in a compact subset of  $E$ .
2.  $T: X \mapsto E$  is expansive and onto.

Then there exists a point  $x^* \in X$  such that
$$Sx^* + Tx^* = x^*.$$

**Definition 2.2** *Let  $(X, d)$  be a metric space and  $M$  be a subset of  $X$ . The mapping  $T: M \mapsto X$  is said to be expansive if there exists a constant  $h > 1$  such that*
$$d(Tx, Ty) \geq h d(x, y)$$
for any  $x, y \in M$ .

## 3   Proof of the Main Result

Let  $\varepsilon \in \left(0, \frac{1}{10}\right)$ .

For  $u \in \mathcal{C}^1([0,1], \mathcal{C}_0^1([0,1]))$ , we define the operators
$$T_{11}u(t,x) = (1+\varepsilon)u(t,x),$$

$$S_{11}u(t,x) = -\varepsilon u(t,x)$$

$$+\varepsilon \left( \int_0^x u(t,z)dz - \int_0^x \phi(z)dz + \frac{1}{2}\int_0^t (u(\tau,x))^2 d\tau \right),$$
$t \in [0,1], \quad x \in [0,1]$ . Consider the IVP
$$u_t + uu_x = 0, \quad (t,x) \in (0,1] \times [0,1],$$

$$u(0,x) = \phi(x), \quad x \in [0,1]. \qquad (3.1)$$

Note that every  $u \in \mathcal{C}^1([0,1], \mathcal{C}_0^1([0,1]))$  which is a fixed point of  $T_{11} + S_{11}$  is a solution of the IVP (3.1). Really, we have

$$u(t,x) = (1+\varepsilon)u(t,x) - \varepsilon u(t,x)$$

$$+ \varepsilon \left( \int_0^x u(t,z)dz - \int_0^x \phi(z)dz + \frac{1}{2}\int_0^t (u(\tau,x))^2 d\tau \right),$$

$(t,x) \in [0,1] \times [0,1]$, whereupon

$$\int_0^x u(t,z)dz - \int_0^x \phi(z)dz + \frac{1}{2}\int_0^t (u(\tau,x))^2 d\tau = 0, \quad (t,x) \in [0,1] \times [0,1]. \qquad (3.2)$$

We differentiate the last equation with respect to $x$ and we get

$$u(t,x) - \phi(x) + \int_0^t u(\tau,x)u_x(\tau,x)d\tau = 0, \quad (t,x) \in [0,1] \times [0,1].$$

Now we differentiate the last equation with respect to $t$ and we obtain

$$u_t(t,x) + u(t,x)u_x(t,x) = 0, \quad (t,x) \in [0,1] \times [0,1].$$

We put $t = 0$ in (3.2). Then

$$\int_0^x u(0,z)dz - \int_0^x \phi(z)dz = 0, \quad x \in [0,1],$$

which we differentiate with respect to $x$ and we get

$$u(0,x) = \phi(x), \quad x \in [0,1].$$

Let $\tilde{\tilde{\tilde{X}}}^1$ be the set of all equicontinuous families of the space $\mathcal{C}^1([0,1], \mathcal{C}_0^1([0,1]))$ with respect to the norm

$$\| f \| == \max\{ \max_{(t,x)\in[0,1]\times[0,1]} |u(t,x)|, \max_{(t,x)\in[0,1]\times[0,1]} |u_t(t,x)|$$

$$\max_{(t,x)\in[0,1]\times[0,1]} |u_x(t,x)|\},$$

$\tilde{X}^1 = \tilde{\tilde{\tilde{X}}}^1 \cup \{\phi\}$, $\tilde{X}^1 = \overline{\tilde{\tilde{\tilde{X}}}^1}$, i.e., $\tilde{X}^1$ is the completion of $\tilde{\tilde{\tilde{X}}}^1$,

$$X^1 = \{f \in \tilde{X}^1 : \| f \| \leq \tfrac{1}{2}\},$$

$$Y^1 = \{f \in \tilde{X}^1 : \| f \| \leq \tfrac{1+\varepsilon}{2}\}.$$

Note that $X^1 \subset Y^1$ and $Y^1$ is a compact subset of $\mathcal{C}^1([0,1], \mathcal{C}_0^1([0,1]))$. For any $x,y \in X^1$, we have

$$\| T_{11}x - T_{11}y \| = (1+\varepsilon) \| x - y \|.$$

Consequently $T_{11} : X^1 \to \mathcal{C}^1([0,1], \mathcal{C}_0^1([0,1]))$ is a linear operator which is expansive with a constant $1+\varepsilon$. Let $v \in Y^1$. Then $u = \frac{v}{1+\varepsilon} \in X^1$ and $T_{11}u = v$. Therefore $T_{11} : X^1 \to Y^1$ is onto. Now we will prove that $S_{11} : X^1 \to X^1$. Let $u \in X^1$. We have

$$|S_{11}u(t,x)| \leq \varepsilon(|u(t,x)| + \int_0^x |u(t,z)|dz + \int_0^x |\phi(z)|dz$$

$$+ \frac{1}{2}\int_0^t |u(\tau,x)|^2 d\tau)$$

$$\leq \frac{13}{8}\varepsilon$$

$$\leq \frac{1}{2}, \quad (t,x) \in [0,1] \times [0,1],$$

and

$$\frac{\partial}{\partial t}(S_{11}u(t,x)) = \varepsilon\left(-u_t(t,x) + \int_0^x u_t(t,z)dz + \frac{1}{2}(u(t,x))^2\right),$$

$$|\frac{\partial}{\partial t}(S_{11}u(t,x))| \leq \varepsilon \left(|u_t(t,x)| + \int_0^x |u_t(t,z)|dz + \frac{1}{2}(u(t,x))^2\right)$$

$$\leq \frac{9}{8}\varepsilon$$

$$\leq \frac{1}{2},$$

and

$$\frac{\partial}{\partial x}(S_{11}u(t,x)) = \varepsilon \left(-u_x(t,x) + u(t,x) - \phi(x) + \int_0^t u(\tau,x)u_x(\tau,x)d\tau\right),$$

$$\left|\frac{\partial}{\partial x}(S_{11}u(t,x))\right| \leq \varepsilon \left(|u_x(t,x)| + |u(t,x)| + |\phi(x)| + \int_0^t |u(\tau,x)u_x(\tau,x)|d\tau\right)$$

$$\leq \frac{7}{4}\varepsilon$$

$$\leq \frac{1}{2}.$$

Consequently $S_{11}: X^1 \to X^1$. Hence and Theorem 2.1, it follows that the IVP (3.1) has a solution $u^{11} \in \mathcal{C}^1([0,1], \mathcal{C}_0^1([0,1]))$. Now, for $u \in \mathcal{C}^1([0,1], \mathcal{C}_0^1([1,2]))$, we define the operators

$$T_{12}u(t,x) = (1+\varepsilon_1)u(t,x),$$

$$S_{12}u(t,x) = -\varepsilon_1 u(t,x)$$

$$+\varepsilon \left(\int_1^x u(t,z)dz - \int_1^x \phi(z)dz + \frac{1}{2}\int_0^t ((u(\tau,x))^2 - (u^{11}(\tau,1))^2)d\tau\right),$$

$t \in [0,1]$, $x \in [1,2]$, where $\varepsilon_1 \in \left(0, \frac{1}{10^2}\right)$. Consider the IVP

$$u_t + uu_x = 0, \quad (t,x) \in [0,1] \times [1,2],$$

$$u(0,x) = \phi(x), \quad x \in [1,2].$$
(3.3)

Note that every $u \in \mathcal{C}^1([0,1], \mathcal{C}_0^1([1,2]))$ which is a fixed point of $T_{12} + S_{12}$ is a solution of the IVP (3.3). Let $\tilde{\tilde{X}}^2$ be the set of all equicontinuous families of the space $\mathcal{C}^1([0,1], \mathcal{C}_0^1([1,2]))$ with respect to the norm

$$\| f \| == \max\{\max_{(t,x)\in[0,1]\times[1,2]}|u(t,x)|, \max_{(t,x)\in[0,1]\times[1,2]}|u_t(t,x)|$$

$$\max_{(t,x)\in[0,1]\times[1,2]}|u_x(t,x)|\},$$

$\tilde{X}^2 = \tilde{\tilde{X}}^2 \cup \{\phi, u^{11}\}$, $\bar{X}^2 = \overline{\tilde{X}^2}$, i.e., $\bar{X}^2$ is the completion of $\tilde{X}^2$,

$$X^2 = \{f \in \bar{X}^2 : \| f \| \leq \frac{1}{2^2}\},$$

$$Y^2 = \{f \in \bar{X}^2 : \| f \| \leq \frac{1+\varepsilon_1}{2^2}\}.$$

Note that $X^2 \subset Y^2$ and $Y^2$ is a compact subset of $\mathcal{C}^1([0,1], \mathcal{C}_0^1([1,2]))$. As in above we get a solution $u^{12} \in \mathcal{C}^1([0,1], \mathcal{C}_0^1([1,2]))$ of the IVP (3.3). Note that

$$\int_1^x u^{12}(t,z)dz - \int_1^x \phi(z)dz + \frac{1}{2}\int_0^t ((u^{12}(\tau,x))^2 - (u^{11}(\tau,1))^2)d\tau = 0,$$

$(t,x) \in [0,1] \times [0,1]$. We have
$$u^{12}(t,1) = u^{11}(t,1) = 0, \quad u_x^{12}(t,1) = u_x^{11}(t,1) = 0,$$
whereupon
$$u_t^{11}(t,1) = u_t^{12}(t,1).$$
Thus, the function
$$v(t,x) = \begin{cases} u^{11}(t,x) & (t,x) \in [0,1] \times [0,1] \\ u^{12}(t,x) & (t,x) \in [0,1] \times [1,2] \end{cases}$$
is a $\mathcal{C}^1([0,1], \mathcal{C}_0^1([0,2]))$-solution of the IVP
$$u_t + uu_x = 0 \quad (t,x) \in [0,1] \times [0,2]$$
$$u(0,x) = \phi(x), \quad x \in [0,2].$$
Next, for $u \in \mathcal{C}^1([0,1], \mathcal{C}_0^1([2,3]))$, we define the operators
$$T_{13}u(t,x) = (1+\varepsilon_2)u(t,x),$$
$$S_{13}u(t,x) = -\varepsilon_2 u(t,x)$$
$$+\varepsilon_2(\int_2^x u(t,z)dz - \int_2^x \phi(z)dz$$
$$+\frac{1}{2}\int_0^t ((u(\tau,x))^2 - (u^{12}(\tau,2))^2)d\tau),$$
$t \in [0,1], \quad x \in [2,3]$, where $\varepsilon \in \left(0, \frac{1}{10^4}\right)$. Consider the IVP
$$u_t + uu_x = 0, \quad (t,x) \in [0,1] \times [2,3],$$
$$u(0,x) = \phi(x), \quad x \in [2,3]. \tag{3.4}$$
Note that every $u \in \mathcal{C}^1([0,1], \mathcal{C}_0^1([2,3]))$ which is a fixed point of $T_{13} + S_{13}$ is a solution of the IVP (3.4). Let $\tilde{\tilde{X}}^3$ be the set of all equicontinuous families of the space $\mathcal{C}^1([0,1], \mathcal{C}_0^1([2,3]))$ with respect to the norm
$$\| f \| == \max\{\max_{(t,x) \in [0,1] \times [2,3]} |u(t,x)|, \max_{(t,x) \in [0,1] \times [2,3]} |u_t(t,x)|$$
$$\max_{(t,x) \in [0,1] \times [2,3]} |u_x(t,x)|\},$$
$$\tilde{X}^3 = \tilde{\tilde{X}}^3 \cup \{\phi, u^{12}\}, \quad \bar{X}^3 = \overline{\tilde{\tilde{X}}^3}, \text{ i.e., } \tilde{X}^3 \text{ is the completion of } \tilde{\tilde{X}}^3,$$
$$X^3 = \{f \in \bar{X}^2 : \| f \| \le \frac{1}{2^3}\},$$
$$Y^2 = \{f \in \bar{X}^2 : \| f \| \le \frac{1+\varepsilon_2}{2^3}\}.$$
Note that $X^3$ is a compact subset of $Y^3$ and $Y^3$ is a compact subset of $\mathcal{C}^1([0,1], \mathcal{C}_0^1([2,3]))$. As in above we get a solution $u^{13} \in \mathcal{C}^1([0,1], \mathcal{C}_0^1([2,3]))$ of the IVP (3.4). Note that
$$\int_2^x u^{13}(t,z)dz - \int_2^x \phi(z)dz + \frac{1}{2}\int_0^t ((u^{13}(\tau,x))^2 - (u^{12}(\tau,2))^2)d\tau = 0,$$
$(t,x) \in [0,1] \times [0,1]$. We have
$$u^{12}(t,2) = u^{13}(t,2) = 0, \quad u_x^{13}(t,2) = u_x^{12}(t,2) = 0,$$

whereupon
$$u_t^{13}(t,2) = u_t^{12}(t,2).$$
Thus, the function
$$v(t,x) = \begin{cases} u^{11}(t,x) & (t,x) \in [0,1] \times [0,1] \\ u^{12}(t,x) & (t,x) \in [0,1] \times [1,2] \\ u^{13}(t,x) & (t,x) \in [0,1] \times [2,3] \end{cases}$$
is a $\mathcal{C}^1([0,1], \mathcal{C}_0^1([0,3]))$-solution of the IVP
$$u_t + uu_x = 0 \quad (t,x) \in [0,1] \times [0,3]$$
$$u(0,x) = \phi(x), \quad x \in [0,3].$$
Continuing in this manner, we get a $\mathcal{C}^1([0,1], \mathcal{C}^1(\mathbb{R}))$-solution $u^1$ of the IVP
$$u_t + uu_x = 0, \quad (t,x) \in [0,1] \times \mathbb{R},$$
$$u(0,x) = \phi(x), \quad x \in \mathbb{R},$$
for which we have
$$|u^1|, \quad |u_t^1|, \quad |u_x^1| \leq \frac{1}{2^{|j|+1}} \quad on \quad [0,1] \times [j, j+1], \quad j \in \mathbb{Z}.$$

Next, for $u \in \mathcal{C}^1([1,2], \mathcal{C}_0^1([0,1]))$, we define the operators
$$T_{21}u(t,x) = (1+\varepsilon)u(t,x),$$
$$S_{21}u(t,x) = -\varepsilon u(t,x)$$
$$+\varepsilon\left(\int_0^x u(t,z)dz - \int_0^x u^1(1,z)dz + \frac{1}{2}\int_0^t (u(\tau,x))^2 d\tau\right),$$
$t \in [0,1]$, $x \in [0,1]$, where $\varepsilon$ is as in Step 1. Consider the IVP
$$u_t + uu_x = 0, \quad (t,x) \in [1,2] \times [0,1],$$
$$u(1,x) = \phi(x), \quad x \in [0,1]. \tag{3.5}$$
Note that every $u \in \mathcal{C}^1([1,2], \mathcal{C}_0^1([0,1]))$ which is a fixed point of $T_{21} + S_{21}$ is a solution of the IVP (3.5). Let $\tilde{\tilde{X}}^{21}$ be the set of all equicontinuous families of the space $\mathcal{C}^1([1,2], \mathcal{C}_0^1([0,1]))$ with respect to the norm
$$\| f \| = \max\{ \max_{(t,x)\in[1,2]\times[0,1]} |u(t,x)|, \quad \max_{(t,x)\in[1,2]\times[0,1]} |u_t(t,x)|$$
$$\max_{(t,x)\in[1,2]\times[0,1]} |u_x(t,x)|\},$$
$\tilde{X}^{21} = \tilde{\tilde{X}}^{21} \cup \{u^1\}$, $\bar{\tilde{X}}^{21} = \overline{\tilde{\tilde{X}}^{21}}$, i.e., $\bar{\tilde{X}}^{21}$ is the completion of $\tilde{\tilde{X}}^{21}$,
$$X^{21} = \{f \in \bar{\tilde{X}}^{21}: \| f \| \leq \frac{1}{2}\},$$
$$Y^{21} = \{f \in \tilde{X}^{21}: \| f \| \leq \frac{1+\varepsilon}{2}\}.$$
Note that $X^{21} \subset Y^{21}$ and $Y^{21}$ is a compact subset of $\mathcal{C}^1([1,2], \mathcal{C}_0^1([0,1]))$. As in above we get a solution $u^{21} \in \mathcal{C}^1([1,2], \mathcal{C}_0^1([0,1]))$ of the IVP (3.5). Note that

$$\int_0^x u^{21}(t,z)dz - \int_0^x u^1(1,z)dz + \frac{1}{2}\int_1^t (u^{21}(\tau,x))^2 d\tau = 0, \quad (t,x) \in [1,2] \times [0,1].$$

We put $t = 1$ in the last equation, then we differentiate with respect to $x$ and we get

$$u^{21}(1,x) = u^1(1,x), \quad x \in [0,1].$$

Hence,

$$u_x^{21}(1,x) = u_x^1(1,x), \quad x \in [0,1],$$

and then

$$u_t^{21}(1,x) = u_t^1(1,x), \quad x \in [0,1].$$

And so on, we get a $\mathcal{C}^1([1,2], \mathcal{C}_0^1(\mathbb{R}))$-solution $u^2$ of the IVP

$$u_t + u u_x = 0, \quad (t,x) \in [1,2] \times \mathbb{R},$$

$$u(0,x) = u^1(1,x), \quad x \in \mathbb{R},$$

for which we have

$$|u^2|, \quad |u_t^2|, \quad |u_x^2| \leq \frac{1}{\sqrt{2}^{|j|+1}} \quad on \quad [1,2] \times [j, j+1], \quad j \in \mathbb{Z}.$$

And so on. This completes the proof.

## 4    Conclusions

In this paper we investigate for existence of classical solutions of a class of Burgers equations. We propose a new approach for investigation for existence of classical solutions. The proposed approach can be used for other classes partial differential equations.